\documentclass[11pt]{article}
\usepackage{amssymb, amsmath, url, graphicx, setspace, geometry}
\usepackage{calrsfs}
\usepackage{wasysym}
\usepackage{xcolor}

\def\3{\subset }
\def\4{\subseteq }
\def\<{\left<}
\def\>{\right>}

\def\bit{\begin{itemize}}
\def\eit{\end{itemize}}
\def\3{\subset }
\def\4{\subseteq }

\def\0{\leqno}

\def\barr{\begin{array}}
\def\earr{\end{array}}

\def\Z{{\rlap{$\kern2pt{\rm Z}$}{\rm Z}\,}}

\title{\bf On the average order of a finite group}
\author{Mihai-Silviu Lazorec and Marius T\u arn\u auceanu}

\begin{document}
\maketitle

\begin{abstract}
Let $o(G)$ be the average order of a finite group $G$. We show that if $o(G)<c$, where $c\in \lbrace \frac{13}{6}, \frac{11}{4}\rbrace$, then $G$ is an elementary abelian 2-group or a solvable group, respectively. Also, we prove that the set containing the average orders of all finite groups is not dense in $[a, \infty)$, for all $a\in [0, \frac{13}{6}]$. We also outline some results related to the integer values of the average order. Since group element orders is  a popular research topic, we pose some open problems concerning the average order of a finite group throughout the paper.
\end{abstract}



\section{Introduction}

Let $G$ be a finite group and let $n\geq 2$ be an integer. Denote by $o(x)$ and $C_n$ the order of an element $x$ in $G$ and the cyclic group of order $n$, respectively. The average order of $G$ is a quantity denoted by $o(G)$ that is given by
$$o(G)=\frac{\psi(G)}{|G|},$$
where $\psi(G)=\sum\limits_{x\in G}o(x)$ is the sum of element orders of $G$. In [9], A. Jaikin-Zapirain uses the average order concept while determining a lower bound for the number of conjugacy classes of a finite $p$-group/nilpotent group $G$. In the same paper, the author suggests the following question: \textit{``Let $G$ be a finite ($p$-)group and $N$ a normal (abelian) subgroup of $G$. Is it true that $o(G)\geq o(N)^{\frac{1}{2}}?$"}. Recently, E.I. Khukhro, A. Moret\' o and M. Zarrin proved the following result (see Theorem 1.2 of \cite{11}):\\

\textbf{Theorem 1.1.} \textit{Let $c>0$ be a real number and let $p\geq \frac{3}{c}$ be a prime. Then there exists a finite $p$-group $G$ with a normal abelian subgroup $N$ such that $o(G)<o(N)^c$.}\\

Note that Theorem 1.1 provides a negative answer to Jaikin-Zapirain's question even if we replace the exponent $\frac{1}{2}$ with any positive real number $c$. 

We can not overlook the fact that the average order $o(G)$ of a group $G$ is strongly related to the sum of element orders $\psi(G)$. During the last years, the latter concept was frequently studied and became a fruitful research topic. Some relevant results involving the sum of element orders are summed up into the following two theorems:\\

\textbf{Theorem 1.2.} \textit{Let $G$ be a finite group and let $\psi'(G)=\frac{\psi(G)}{\psi(C_{|G|})}$.
\begin{itemize}
\item[a)] If $\psi'(G)>\frac{7}{11}=\psi'(C_2\times C_2)$, then $G$ is cyclic; (see Theorem 1 of \cite{7})
\item[b)] If $\psi'(G)>\frac{13}{21}=\psi'(S_3)$, then $G$ is nilpotent; (see Theorem 1.1 of \cite{15})
\item[c)] If $\psi'(G)>\frac{31}{77}=\psi'(A_4)$, then $G$ is supersolvable; (see the proof of [T]-Conjecture in \cite{4})
\item[d)] If $\psi'(G)>\frac{211}{1617}=\psi'(A_5)$, then $G$ is solvable. (see the main result of \cite{3})
\end{itemize}}

\textbf{Theorem 1.3.} \textit{(see Theorem 1.1 of \cite{14}) Let $G$ be a finite group and let $\psi''(G)=\frac{\psi(G)}{|G|^2}$.
\begin{itemize}
\item[a)] If $\psi''(G)>\frac{7}{16}=\psi''(C_2\times C_2)$, then $G$ is cyclic;
\item[b)] If $\psi''(G)>\frac{27}{64}=\psi''(Q_8)$, then $G$ is abelian;
\item[c)] If $\psi''(G)>\frac{13}{36}=\psi''(S_3)$, then $G$ is nilpotent;
\item[d)] If $\psi''(G)>\frac{31}{144}=\psi''(A_4)$, then $G$ is supersolvable;
\item[e)] If $\psi''(G)>\frac{211}{3600}=\psi''(A_5)$, then $G$ is solvable.
\end{itemize}}

Note that the lower bounds from Theorems 1.2 and 1.3 are the best possible ones and all items actually constitute some criteria for the cyclicity, commutativity, nilpotency, (super)solvability of a finite group. Other interesting aspects on the sum of element orders concern its maximum and minimum values. Hence, we would like to recall the main result of \cite{2} which is due to H. Amiri, S.M. Jafarian Amiri and I.M. Isaacs. It states that given a finite group $G$ of order $n$, we have
$$\psi(G)\leq \psi(C_n) \text{ and } \psi(G)=\psi(C_n)\Longleftrightarrow G\cong C_n.$$
In other words, the maximum value of $\psi$ among finite groups of order $n$ is attained by $C_n$. If we let $G$ to be a non-cyclic group of order $n$, it is known that 
$$\psi(G)\leq \frac{7}{11}\psi(C_n),$$
this result being proved by M. Herzog, P. Longobardi and M. Maj in \cite{8}. In what concerns the minimum value of $\psi$, the investigation is still open. A complete result is known if one works with nilpotent groups. More exactly, in \cite{1}, H. Amiri and S.M. Jafarian Amiri show that among all nilpotent groups of order $n$, the minimum value of $\psi$ is attained by a group whose Sylow subgroups have prime exponents. 
 
Our main aim is to establish some criteria, as the ones included in Theorems 2 and 3, using the average order of a finite group $G$. The elements of order 2 will play a significant role in our study. Hence, let
$$n_2(G)=|\lbrace x\in G \ | \ o(x)=2\rbrace|$$
be the number of elements of order 2 in $G$. We recall some results which connect this number with the nature/structure of $G$.\\

\textbf{Theorem 1.4.} \textit{Let $G$ be a finite group and let $n_2(G)$ be the number of elements of order 2 in $G$.
\begin{itemize}
\item[a)] If $n_2(G)\geq \frac{3}{4}|G|$, then $G$ is an elementary abelian 2-group;
\item[b)] If $n_2(G)>\frac{1}{4}|G|-1$, then $G$ is solvable or there exists $m\in\mathbb{N}$ such that $G\cong C_2^m\times A_5$.
\end{itemize}}

Item $a)$ of Theorem 1.4 is a well-known result and a proof can be found in the preprint \cite{5} (see Theorem 4.1). For item $b)$, the reader may consult Theorem 2.44 of \cite{6} or \cite{12}. The same results may be rewritten in terms of the number of solutions of the equation $x^2=1$, where $x\in G$. 

We end this section with a brief description of our main results. We prove the following theorem:\\

\textbf{Theorem 1.5.} \footnote{Meanwhile, the upper bound outlined in item $b)$ of Theorem 1.5 was improved. More exactly, in \cite{9}, it was proved that given a finite group $G$, if $o(G)<\frac{211}{60}=o(A_5)$, then $G$ is solvable. This result was also conjectured in \cite{11}.} \textit{Let $G$ be a finite group.
\begin{itemize}
\item[a)] If $o(G)<\frac{13}{6}=o(S_3)$, then $G$ is an elementary abelian 2-group;
\item[b)] If $o(G)<\frac{11}{4}=o(C_4)$, then $G$ is solvable.
\end{itemize}}

Besides this result, we approach the problem of determining the integer values of the average order of a finite group. More exactly, we show that there are no finite groups $G$ such that $o(G)\in \lbrace 2,3\rbrace$. We also prove that there are infinitely many values located to the right of and close to $\frac{13}{6}$ which are not attained by the average order of a finite group. As a consequence, we show that the set containing the average orders of all finite groups is not dense in $[\frac{13}{6}, \infty)$. The corresponding proofs, other results and some open problems are outlined in the following section.

\section{Main results}

We start by outlining two auxiliary results along with their proofs. Both are needed in order to prove the criteria included in Theorem 1.5.\\

\textbf{Lemma 2.1.} \textit{Let $c>0$ be a real number, let $G$ be a finite group of order $n$, $1=d_1<d_2<\ldots <d_r$ be the element orders of $G$, where $r\geq 3$, and let $n_{d_i}=|\lbrace x\in G \ | \ o(x)=d_i\rbrace|$, for all $i\in\lbrace 1,2,\ldots, r \rbrace$. If $o(G)<c$, then
\begin{align}\label{r1}
n_{d_2}>\frac{d_3-c}{d_3-d_2}n-\frac{d_3-1}{d_3-d_2}.
\end{align}}

\textbf{Proof.} We have $n=1+n_{d_2}+n_{d_3}+\ldots +n_{d_r}$, so we deduce that
\begin{align*}
\psi(G)&=1+d_2n_{d_2}+d_3n_{d_3}+\ldots +d_rn_{d_r}\\ 
&\geq 1+d_2n_{d_2}+d_3(n_{d_3}+\ldots n_{d_r})\\
&=1+d_2n_{d_2}+d_3(n-1-n_{d_2}).
\end{align*}
Since $o(G)<c$, it follows that $\psi(G)<cn$. Therefore,
$$cn>1+d_2n_{d_2}+d_3(n-1-n_{d_2}).$$
The last inequality is equivalent to (\ref{r1}), so our proof is complete.
\hfill\rule{1,5mm}{1,5mm}\\ 

\textbf{Lemma 2.2.} \textit{Let $G$ be a finite group, let $p$ be a prime number and let $m\in \mathbb{N}$. Then
\begin{align}\label{r2}
o(C_p^m\times G)=o(G)+\frac{o(C_p^m)-1}{|G|}\sum\limits_{\substack{ x\in G \\ p \ \nmid \ o(x)}}o(x).
\end{align}}

\textbf{Proof.} We denote the least common multiple of $o(a)$ and $o(x)$ by $[o(a), o(x)]$, for all $(a, x)\in C_p^m\times G$. We have
\begin{align*}
\psi(C_p^m\times G)&=\sum\limits_{(a,x)\in C_p^m\times G}o(a,x)=\sum\limits_{a\in C_p^m}\sum\limits_{x\in G}[o(a), o(x)]\\
&=\sum\limits_{a\in C_p^m}\bigg(\sum\limits_{\substack{ x\in G \\ p \ | \ o(x)}}[o(a), o(x)]+\sum\limits_{\substack{ x\in G \\ p \ \nmid \ o(x)}}[o(a), o(x)]\bigg)\\
&=\sum\limits_{a\in C_p^m}\bigg(\sum\limits_{\substack{ x\in G \\ p \ | \ o(x)}}o(x)+\sum\limits_{\substack{ x\in G \\ p \ \nmid \ o(x)}}o(a)\cdot o(x)\bigg)\\
&=p^m\sum\limits_{\substack{ x\in G \\ p \ | \ o(x)}}o(x)+\bigg(\sum\limits_{a\in C_p^m}o(a) \bigg)\cdot \bigg(\sum\limits_{\substack{ x\in G \\ p \ \nmid \ o(x)}}o(x) \bigg)=\\
&=p^m\bigg(\psi(G)-\sum\limits_{\substack{ x\in G \\ p \ \nmid \ o(x)}}o(x)\bigg)+\psi(C_p^m)\sum\limits_{\substack{ x\in G \\ p \ \nmid \ o(x)}}o(x)=\\
&=p^m\psi(G)+(\psi(C_p^m)-p^m)\sum\limits_{\substack{ x\in G \\ p \ \nmid \ o(x)}}o(x),
\end{align*}
and the conclusion easily follows.
\hfill\rule{1,5mm}{1,5mm}\\ 

We proceed by justifying the results that are summed up into Theorem 1.5.\\

\textbf{Proof of Theorem 1.5.} Let $G$ be a finite group of order $n$. We use the notations of Lemma 2.1 in our reasoning.

\textit{a)} If $r=2$, then the set of element orders of $G$ is $\pi_e(G)=\lbrace 1,p\rbrace$ and $n=p^k$, where $p$ is a prime number and $k\geq 1$. For $p=2$, it follows that $G$ is an elementary abelian 2-group. For $p\geq 3$, we have
$$o(G)=\frac{1+p(p^k-1)}{p^k}=p-\frac{p-1}{p^k}\geq p-\frac{p-1}{p}\geq 3-\frac{2}{3}>\frac{13}{6},$$
a contradiction. 

Therefore, we may assume that $r\geq 3$ and apply Lemma 2.1. If $d_2\geq 3$, then, by following the steps above, we get $o(G)\geq 3-\frac{2}{3}>\frac{13}{6}$, a contradiction. If $d_2=2$, then Lemma 2.1 leads to
$$n_2>\frac{d_3-\frac{13}{6}}{d_3-2}n-\frac{d_3-1}{d_3-2}\geq \frac{3-\frac{13}{6}}{3-2}n-\frac{3-1}{3-2}=\frac{5}{6}n-2.$$
For $n\geq 24$, we have $\frac{5}{6}n-2\geq \frac{3}{4}n$, so $n_2>\frac{3}{4}n$. According to Theorem 1.4, \textit{a)}, we deduce that $G$ is an elementary abelian 2-group. To investigate the case $n\leq 23$, we use GAP \cite{16} to check that the quantity
$$a_n=\min\lbrace o(G) \ | \ |G|=n\rbrace$$
attains the following values:
$$a_1=1, a_2=\frac{3}{2}, a_3=\frac{7}{3}, a_4=\frac{7}{4}, a_5=\frac{21}{5}, a_6=\frac{13}{6}, a_7=\frac{43}{7}, a_8=\frac{15}{8}, a_9=\frac{25}{9}, a_{10}=\frac{31}{10}$$
$$a_{11}=\frac{111}{11}, a_{12}=\frac{31}{12}, a_{13}=\frac{157}{13}, a_{14}=\frac{57}{14}, a_{15}=\frac{147}{15}, a_{16}=\frac{31}{16}, a_{17}=\frac{273}{17}, a_{18}=\frac{43}{18},$$ 
$$a_{19}=\frac{343}{19}, a_{20}=\frac{71}{20}, a_{21}=\frac{85}{21}, a_{22}=\frac{133}{22}, a_{23}=\frac{507}{23}.$$
Observe that $a_n\geq \frac{13}{6}$ excepting the case $n\in \lbrace 1,2,4,8,16\rbrace$. If $n\in \lbrace 1,2\rbrace$, then $G$ is an elementary abelian 2-group. For $n=4$, we have $o(C_4)=\frac{11}{4}>\frac{13}{6}$. For $n=8$, we get $o(C_8)=\frac{43}{8}>\frac{13}{6}, o(C_2\times C_4)=\frac{23}{8}>\frac{13}{6}, o(Q_8)=\frac{27}{8}>\frac{13}{6}$ and $o(D_8)=\frac{19}{8}>\frac{13}{6}$. Similarly, for $n=16$, we obtain $o(G)>\frac{13}{6}$ for any group $G$ of order 16 such that $G\not\cong C_2^4$. Consequently, the conclusion also holds for $n\in\lbrace 4, 8, 16\rbrace$ and this completes the proof of item \textit{a)}.

Before justifying that criterion \textit{b)} of Theorem 1.5 holds, we would like to mention that the case $n\leq 23$ above may be also solved by proving the following two results: \textit{``If $G$ is a finite $p$-group such that $o(G)<\frac{13}{6}$, then $G$ is an elementary abelian $2$-group."} and \textit{``If $G$ is a finite group such that $|G|\in \lbrace pq, p^2q\rbrace$, where $p$ and $q$ are distinct primes, then $o(G)\geq \frac{13}{6}$ and the equality holds if and only if $G\cong S_3$."} However, we consider that the use of GAP is preferable and more efficient.\\

\textit{b)} If $d_2\geq 3$, then $|G|$ is odd and, consequently, $G$ is solvable. If $r\leq 3$, then $|G|$ is divisible by at most two distinct primes, so $G$ would be solvable once again. Consequently, we may assume that $d_2=2$ and $r\geq 4$. 

In the proof of Lemma 2.1, we showed that
\begin{align}\label{r3}
\psi(G)=1+d_2n_{d_2}+d_3n_{d_3}+\ldots +d_rn_{d_r}\geq 1+d_2n_{d_2}+d_3(n_{d_3}+\ldots n_{d_r}).
\end{align}
Since $r\geq 4$, we have $d_4>d_3$, so an equality can not occur in (\ref{r3}). Therefore,
$$\psi(G)\geq 2+d_2n_{d_2}+d_3(n_{d_3}+\ldots n_{d_r}).$$
If $o(G)<c$, where $c>0$ is a real number, then one obtains the following inequality which is similar to the one given by Lemma 2.1:
$$n_2>\frac{d_3-c}{d_3-d_2}n-\frac{d_3-2}{d_3-d_2}.$$
In our case, for $c=\frac{11}{4}$, we get
$$n_2>\frac{d_3-\frac{11}{4}}{d_3-2}n-\frac{d_3-2}{d_3-2}\geq \frac{3-\frac{11}{4}}{3-2}n-1=\frac{1}{4}n-1.$$
By Theorem 1.4 \textit{b)}, we deduce that $G$ is solvable or there exists $m\in\mathbb{N}$ such that $G\cong C_2^m\times A_5$. To complete our proof, it remains to show that the second case can not occur. Hence, suppose that $G\cong C_2^m\times A_5$. Then, by Lemma 2.2, we obtain
\begin{align*}
o(G)&=o(C_2^m\times A_5)=o(A_5)+\frac{o(C_2^m)-1}{60}\sum\limits_{\substack{ x\in A_5 \\ 2 \ \nmid \ o(x)}}o(x)\\
&=o(A_5)+\frac{181}{60}\cdot\frac{2^{m}-1}{2^m}>o(A_5)=\frac{211}{60}>\frac{11}{4},
\end{align*}
a contradiction. 
\hfill\rule{1,5mm}{1,5mm}\\ 

Notice that item \textit{a)} of Theorem 1.5 may be also viewed as a nilpotency criterion as follows: \textit{``If $o(G)<\frac{13}{6}=o(S_3)$, then $G$ is nilpotent."} Since $S_3$ is non-nilpotent, the ratio $\frac{13}{6}$ is the best upper bound that may appear in such a criterion. Below, we will explain that $\frac{13}{6}$ is also the best upper bound for the original form of item \textit{a)}. 

Further, we outline some consequences of Theorem 1.5. First of all, notice that if $G\cong C_2^m$, where $m\in\mathbb{N}$, we have 
\begin{align}\label{r4}
o(G)=2-\frac{1}{2^m}<2<\frac{13}{6}.
\end{align} 
Hence, by Theorem 1.5, \textit{a)}, we obtain the following characterization of elementary abelian 2-groups.\\

\textbf{Corollary 2.3.} \textit{A finite group $G$ is an elementary abelian 2-group if and only if $o(G)<\frac{13}{6}.$}\\

The inequality (\ref{r4}) shows that all elementary abelian 2-groups have low average orders contained in $[1, 2)$. If we exclude these groups, we can state the following result which outlines some minimum values of the average order on the class of finite $p$-groups. The proof follows some ideas that were already highlighted in this paper, so we omit explicit details.\\

\textbf{Corollary 2.4.} \textit{Let $G$ be a finite $p$-group.
\begin{itemize}
\item[a)] If $p=2$ and $G$ is not elementary abelian, then $o(G)\geq \frac{19}{8}=o(D_8)$;
\item[b)] If $p$ is odd, then $o(G)\geq \frac{p^2-p+1}{p}=o(C_p).$
\end{itemize}}

As a consequence of inequality (\ref{r4}) and Corollary 2.3, we can state that the upper bound $\frac{13}{6}$ in Theorem 1.5, \textit{a)}, is the best possible one. Also, we deduce a first result concerning the  possibility of obtaining specific integer values when we compute the average order of a finite group.\\

\textbf{Corollary 2.5.} \textit{There are no finite groups $G$ such that $o(G)\in [2, \frac{13}{6})$. In particular, there are no finite groups $G$ such that $o(G)=2$.}\\

Another way to justify that $o(G)\ne 2$, for any finite group $G$, is by noticing that the sum of element orders of a finite group is always an odd positive integer. We use this remark to show the non-existence of finite groups whose average orders are equal to 3.\\

\textbf{Proposition 2.6.} \textit{There are no finite groups $G$ such that $o(G)=3$.}\\

\textbf{Proof.} Suppose that there is a finite group $G$ such that $o(G)=3$. Since $\psi(G)$ is odd and $o(G)$ is an integer, it follows that $n=|G|$ is also odd. If $exp(G)=3$, then
$$o(G)=\frac{1+3(n-1)}{n}=3-\frac{2}{n}<3.$$
Otherwise, $G$ has at least 4 elements of order greater than or equal to 5, so
$$o(G)\geq \frac{1+4\cdot 5+3\cdot (n-5)}{n}=3+\frac{6}{n}>3.$$
Hence, in both cases we arrive at a contradiction. Consequently, our proof is complete. 
\hfill\rule{1,5mm}{1,5mm}\\

Taking into account the previous two results, an interesting problem would be to determine the integer values contained in $Im(o)$, where  
$$Im(o)=\lbrace o(G) \ | \ G\in\mathcal{G}\rbrace$$
and $\mathcal{G}$ is the class of all finite groups. We note that such values exist. Besides the trivial case, GAP yields the following examples of groups whose orders are at most 5000 and whose average orders are integers:
\begin{itemize}
\item[--] $G_1\cong C_5\times (C_7\rtimes C_3) \ (SmallGroup(105,1)), o(G_1)=\frac{1785}{105}=17$;
\item[--] $G_2\cong C_{17}\times (C_7\rtimes C_3) \ (SmallGroup(357,1)), o(G_2)=\frac{23205}{357}=65;$
\item[--] $G_3\cong C_{85}\times (C_7\rtimes C_3) \ (SmallGroup(1785,1)), o(G_3)=\frac{487305}{1785}=273$;
\item[--] $G_4\cong C_{299}\times C_{13} \ (SmallGroup(3887,2)), o(G_4)=\frac{13446147}{3887}=285$;
\item[--] $G_5\cong C_{35}\times (C_{43}\rtimes C_3) \ (SmallGroup(4515,3)), o(G_5)=\frac{1864695}{4515}=413$;
\item[--] $G_6\cong C_{221}\times (C_7\rtimes C_3) \ (SmallGroup(4641,3)), o(G_6)=\frac{3643185}{4641}=785.$
\end{itemize}
Hence, we consider that it is natural to pose the following question.\\

\textbf{Question 2.7.} \textit{Which are the integer values contained in $Im(o)$? Moreover, which are the finite groups $G$ having an element of average order, i.e. there exists $x\in G$ such that $o(G)=o(x)$?}\\

The final part of this section concerns the study of the density of $Im(o)$ in $[\frac{13}{6}, \infty)$. First of all, we justify the choice of the interval $[\frac{13}{6}, \infty)$. If one would opt for an interval such as $I=[a, \infty)$, where $a<\frac{13}{6}$, then Corollary 2.5 assures that $Im(o)$ is not dense in $I$. Consequently, we are interested in intervals such as $[a, \infty)$, with $a\geq \frac{13}{6}$. Secondly, by Lemma 2.1 of \cite{1}, we know that if $G_1$ and $G_2$ are two finite groups of coprime orders, then $\psi(G_1\times G_2)=\psi(G_1)\psi(G_2)$. This property is obviously inherited by the average order. We use this multiplicativity property in the proof of the following result which states that any integer $n\geq 2$ is an adherent point of $Im(o)$.\\

\textbf{Proposition 2.8.} \textit{Let $n\geq 2$ be an integer. Then there is a sequence $(G_m)_{m\geq 1}$ of finite groups such that $\displaystyle\lim_{m\to\infty} o(G_m)=n.$}\\

\textbf{Proof.} Let $n=p_1^{n_1}p_2^{n_2}\ldots p_k^{n_k}$ be the prime factorization of $n$. For $i\in\lbrace 1, 2,\ldots, k\rbrace$, we consider the finite group $G_{i,m}=C_{p_i^{n_i}}^{m}$, where $m\geq 1$ is an integer. By Corollary 4.4 of \cite{13}, $G_{i,m}$ has $p_i^{m\alpha}-p_i^{m(\alpha-1)}$ elements of order $p_i^{\alpha}$, for all $\alpha\in\lbrace 1, 2, \ldots, n_i\rbrace$. Hence,
\begin{align*}
\psi(G_{i,m})&=1+p_i(p_i^{m}-1)+p_i^2(p_i^{2m}-p_i^{m})+\ldots +p_i^{n_i}(p_i^{mn_i}-p_i^{m(n_i-1)})\\
&=\frac{p_i^{(m+1)(n_i+1)}-p_i^{(m+1)n_i+1}+p_i-1}{p_i^{m+1}-1},
\end{align*}
and
\begin{align*}
o(G_{i,m})&=\frac{p_i^{(m+1)(n_i+1)}-p_i^{(m+1)n_i+1}+p_i-1}{p_i^{mn_i}(p_i^{m+1}-1)}=p_i^{n_i}\cdot\frac{1-\frac{1}{p_i^{m}}}{1-\frac{1}{p_i^{m+1}}}+\frac{p_i-1}{p_i^{mn_i}(p_i^{m+1}-1)}.
\end{align*}

Finally, let us consider the sequence $(G_m)_{m\geq 1}$, where $G_m=G_{1,m}\times G_{2,m}\times\ldots \times G_{k,m}$. Using the multiplicativity property of the average order, we have
$$\displaystyle\lim_{m\to\infty}o(G_m)=\displaystyle\lim_{m\to\infty}\prod\limits_{i=1}^k o(G_{i,m})=\prod\limits_{i=1}^k \bigg(\displaystyle\lim_{m\to\infty} o(G_{i,m})\bigg)=\prod\limits_{i=1}^k p_i^{n_i}=n,$$
as desired. 
\hfill\rule{1,5mm}{1,5mm}\\

Our following result shows that there is a finite number of average orders which are close to and higher than the ratio $\frac{13}{6}.$\\

\textbf{Proposition 2.9.} \textit{Let $\varepsilon\in (0,\frac{1}{12})$. Then there is a finite number of finite groups $G$ such that $o(G)\in [\frac{13}{6}, \frac{9}{4}-\varepsilon)$.}\\

\textbf{Proof.} We use the notations of Lemma 2.1 in our reasoning. Let $G$ be a finite group of order $n$ such that $o(G)\in [\frac{13}{6}, \frac{9}{4}-\varepsilon)$. 

Suppose that $2 \nmid n$. If $n=1$, then $o(G)=1$, a contradiction. If $n\geq 3$, then $d_2\geq 3$ and $n>\frac{8}{3+4\varepsilon}$. We have 
$$o(G)\geq \frac{1+3(n-1)}{n}=3-\frac{2}{n}\geq \frac{9}{4}-\varepsilon \text{ for } n\geq \frac{8}{3+4\varepsilon},$$
a contradiction.

Therefore, $2|n$. Then $d_2=2$ and we may assume that $r\geq 3$ (otherwise, $G$ would be an elementary abelian 2-group and $o(G)<2$, a contradiction). By Lemma 2.1, we get
$$n_2>\frac{d_3-(\frac{9}{4}-\varepsilon)}{d_3-2}n-\frac{d_3-1}{d_3-2}\geq \frac{3-(\frac{9}{4}-\varepsilon)}{3-2}n-\frac{3-1}{3-2}=\bigg(\frac{3}{4}+\varepsilon\bigg)n-2.$$
Then there is $n_0\in\mathbb{N}$ such that 
$$\bigg(\frac{3}{4}n+\varepsilon\bigg)n-2\geq \frac{3}{4}n, \forall \ n\geq n_0, \text{ so } n_2>\frac{3}{4}n, \forall \ n\geq n_0.$$
If $n\geq n_0$, by Theorem 1.4, \textit{a)}, we deduce that $G$ is an elementary abelian 2-group, so $o(G)<2$, a contradiction. Hence, $o(G)$ may be contained in $[\frac{13}{6}, \frac{9}{4}-\varepsilon)$ only if $|G|=n<n_0$. This means that, up to isomorphism, there is a finite number of groups whose average orders are included in the above interval.
\hfill\rule{1,5mm}{1,5mm}\\

One may restate Proposition 2.9 by saying that there are infinitely many values $v$, located to the right of $\frac{13}{6}$, for which we can not find a finite group $G$ such that $o(G)=v$. This fact provides an answer concerning the density of $Im(o)$ in the interval $[\frac{13}{6}, \infty).$\\

\textbf{Corollary 2.10} \textit{$Im(o)$ is not dense in $[\frac{13}{6}, \infty).$}\\

The investigation over the density of $Im(o)$ can be further continued. Based on our results, it would be natural to ask for an argument which proves/disproves the following conjecture.\\

\textbf{Conjecture 2.11.} \textit{Let $a\geq 0$ be a real number. Then $Im(o)$ is not dense in $[a, \infty)$.}\\

Our results show that this conjecture is valid for $a\in [0,\frac{13}{6}]$. 

\bigskip\noindent {\bf Acknowledgements.} The authors are grateful to the reviewers for their remarks which improve the previous version of the paper. This work was supported by a grant of the "Alexandru Ioan Cuza" University of Iasi, within the Research Grants program, Grant UAIC, code GI-UAIC-2021-01.

\vspace*{3ex}
\small

\begin{minipage}[t]{7cm}
Mihai-Silviu Lazorec \\
Faculty of  Mathematics \\
"Al.I. Cuza" University \\
Ia\c si, Romania \\
e-mail: {\tt silviu.lazorec@uaic.ro}
\end{minipage}
\hspace{3cm}
\begin{minipage}[t]{7cm}
Marius T\u arn\u auceanu \\
Faculty of  Mathematics \\
"Al.I. Cuza" University \\
Ia\c si, Romania \\
e-mail: {\tt tarnauc@uaic.ro}
\end{minipage}

\begin{thebibliography}{100}
\bibitem{1} Amiri, H., \& Jafarian Amiri, S.M. (2011). Sum of element orders on finite groups of the same order. \textit{J. Algebra Appl. 10, no. 2}, 187-190, https://doi.org/10.1142/S0219498811004057.

\bibitem{2} Amiri, H., Jafarian Amiri, S.M., \& Isaacs, I.M. (2009). Sums of element orders in finite groups. \textit{Comm. Algebra 37}, 2978-2980, https://doi.org/10.1080/00927870802502530.

\bibitem{3} Baniasad Azad, M., \& Khosravi, B. (2018). A criterion for solvability of a finite group by the sum of element orders. \textit{J. Algebra 516}, 115-124, https://doi.org/10.1016/j.jalgebra.2018.09.009.

\bibitem{4} Baniasad Azad, M., \& Khosravi, B. (2022). On two conjectures about the sum of element orders. \textit{Canad. Math. Bull. 65 (1)}, 30-38, https://doi.org/10.4153/S0008439521000047.

\bibitem{5} Edmonds, A.L., \& Norwood, Z.B., Finite groups with many involutions. Preprint, arXiv:0911.1154. 

\bibitem{6} Farrokhi, D.G.M. (2015). On the probability that a group satisfies a law: A survey. \textit{Kyoto Univ. Res. Inform. Repos., 1965}, 158-179.

\bibitem{7} Herzog, M., Longobardi, P., \& Maj, M. (2018). An exact upper bound for sums of element orders in non-cyclic finite groups. \textit{J. Pure Appl. Algebra 222, no. 7}, 1628-1642, https://doi.org/10.1016/j.jpaa.2017.07.015.

\bibitem{8} Herzog, M., Longobardi, P., \& Maj, M. (2021) The second maximal groups with respect to the sum of element orders. \textit{J. Pure Appl. Algebra 225, no. 3}, article ID 106531, https://doi.org/10.1016/j.jpaa.2020.106531.

\bibitem{9} Herzog, M., Longobardi, P., \& Maj, M. (2022) Another criterion for solvability of finite groups. \textit{J. Algebra  597}, 1-23, https://doi.org/10.1016/j.jalgebra.2022.01.005.

\bibitem{10} Jaikin-Zapirain, A. (2011). On the number of conjugacy classes of finite nilpotent groups. \textit{Adv. Math. 227}, 1129-1143, doi:10.1016/j.aim.2011.02.021.

\bibitem{11} Khukhro, E.I., Moret\' o, A., \& Zarrin, M. (2021). The average element order and the number of conjugacy classes of finite groups. \textit{J. Algebra 569}, 1-11, https://doi.org/10.1016/j.jalgebra.2020.11.009.

\bibitem{12} Potter, W.M. (1988). Nonsolvable groups with an automorphism inverting many elements. \textit{Arch. Math. 50}, 292-299, https://doi.org/10.1007/BF01190221.


\bibitem{13} T\u arn\u auceanu, M. (2010). An arithmetic method of counting the subgroups of a finite abelian group. \textit{Bull. Math. Soc. Sci. Math. Roumanie (N.S.), 53(101), no. 4}, 373-386.

\bibitem{14} T\u arn\u auceanu, M. (2020) Detecting structural properties of finite groups by the sum of element orders. \textit{Israel J. Math. 238, no. 2}, 629-637, https://doi.org/10.1007/s11856-020-2033-9.

\bibitem{15} T\u arn\u auceanu, M. (2021). A criterion for nilpotency of a finite group by the sum of element orders. \textit{Comm. Algebra 49 (4)}, 1571-1577, https://doi.org/10.1080/00927872.2020.1840575.

\bibitem{16} The GAP Group (2020). \textit{GAP -- Groups, Algorithms, and Programming, Version 4.11.0}, https://www.gap-system.org.
\end{thebibliography}
\end{document}